\input amstex

\magnification1400

\hoffset=0in

\voffset=0in

\noindent

\space

\vskip0pt plus0pt minus0pt

\centerline{\bf Approximation by Lipschitz functions}

\medskip

\centerline{\bf  L. A. Coburn}

\bigskip

\centerline {\bf Abstract}

\medskip

On any metric space, I provide an intrinsic characterization for the uniform closure of the set of all complex-valued Lipschitz functions.  There are applications to function theory on {\it complete Riemannian manifolds} and, in particular, on {\it bounded symmetric domains (BSD)} in ${\text{\bf C}}^{n}$.

\medskip


\smallskip


\bigskip

{\bf 1. Introduction}.  On any metric space $(X, \beta(\cdot,\cdot))$, we say a complex-valued function $f$ is {\it uniformly continuous } if, for arbitrary real $\epsilon > 0$, and $x, y$ in $X$, there is a real $\delta = \delta(\epsilon) > 0$ so that $|f(x) - f(y)| < \epsilon$ whenever $\beta (x, y) < \delta$.  The set of all uniformly continuous functions on $(X, \beta)$ is denoted by $UC(X)$.  The {\it Lipschitz functions} $Lip(X)$ are the subset of $UC(X)$ with the property that, for all $x, y$ in $X$ and f in $Lip(X)$, $|f(x) - f(y)| \leq C \beta (x, y)$, for some positive constant $C = C(f)$.

\medskip

We will be concerned with an intermediate set of functions, {\it the uniform closure  of } $Lip(X)$, which I denote by $Lip_{c}(X)$.  I will show that $Lip_{c}(X)$ consists precisely of those functions $f$ for which, given any $\epsilon > 0$, there is a constant $C= C(\epsilon)$ so that 

$$
|f(x) - f(y)| < \epsilon + C(\epsilon) \beta (x, y),
\tag *
$$
 
 \medskip
 
 \noindent
 for all $x, y$ in $X$.

 \medskip 
 
 As an application of this result, I give a concise proof of the known equivalence $UC(X) \equiv Lip_{c}(X)$ for the special case of complete (connected) Riemannian manifolds $X$, with metric the usual Riemannian distance function induced by the infinitesimal Riemannian metric.  The prototypical complete Riemannian manifold is real n-dimensional space ${\text{\bf R}}^{n}$ and for $x, y$ in ${\text{\bf R}}^{n}$, we have the usual norm $|x|$ and the Riemannian distance function is just $\beta (x, y) = |x - y|$.  This application holds, in particular, for all {\it bounded symmetric domains (BSD)} $\Omega$ in ${\text{\bf C}}^{n}$.  In this case, stronger results (with a more complicated proof) are known [1]: the real-analytic Lipschitz functions are uniformly dense in $UC(\Omega)$.

\medskip

\newpage

\bigskip

\bigskip

For a definitive, fairly recent, treatment of approximation by Lipschitz functions, see [6].  The result characterizing $Lip_{c}(X)$ does not seem to be in the literature.  While evidently not as useful as the notion of ``Lipschitz in the small,'' it still seems to be worth some attention.

 \bigskip

 {\bf 2.  A characterization of} $Lip_{c}(X)$.  For $(X, \beta(\cdot, \cdot))$ any metric space, we recall the extension result due to E. J. McShane [8] for {\it real-valued} Lipschitz functions:
 
 \medskip
 
 \noindent
 {\bf Proposition 1.} For any non-empty subset $S$ of $X$ and any real-valued function $f$ in $Lip(S)$, there is a real-valued function $F$ in $Lip(X)$ with $F|_{S} = f$ and with  $F$ having the same Lipschitz constant as $f$.
 
 \medskip
 
 {\bf Proof.}  Suppose that $|f(s) - f(t)| \leq C \beta(s, t)$ for all $s, t$ in $S$.  For any $x$ in $X$ we define
 
 $$
 F(x)= \inf \{ f(s) + C \beta(x, s): s \in S \}.
 $$
 
 \medskip
 
 \noindent
 To see that $F(x)$ is finite for every $x$ in $X$, fix $s_{0}$ in $S$.  Then we can check that for {\it any} $s$ in $S$
 
 $$
 \align
 f(s) + C\beta(s, x) &\geq f(s_{0}) - C\beta(s, s_{0}) + C\beta(s, x)\\
 &\geq f(s_{0}) - C\beta(x, s_{0}).
 \endalign
 $$
 
 \medskip
 
 \noindent
 For $x$ in $S$, $F(x) \leq f(x)$.  But, for all $s$ in $S$, $f(x) \leq f(s) + C\beta(x, s)$ so $f(x) \leq F(x)$.
 
 \medskip
 
 Finally, we check that $F$ is Lipschitz on $X$.  For $x, y$ in $X$, note that
 
 $$
 \align
 F(x) &= \inf \{ f(s) + C\beta(s, x): s \in S\}\\
 & \leq \inf \{ f(s) + C\beta(s, y) + C\beta (y, x): s \in S\}\\
 &\leq \inf \{ f(s) + C \beta(s, y); s \in S \} + C \beta (y, x) \\
 & \leq F(y) + C \beta(x, y)
 \endalign
 $$
 
 \medskip
 
 \noindent
 so $|F(x) - F(y)| \leq C \beta(x, y)$.
 
 \medskip
 
 \newpage
 
 \bigskip
 
 \bigskip
 
 \noindent
 {\bf Corollary}.  For $S$ any non-empty subset of $X$ and $f$ any {\it complex-valued} function in $Lip(S)$ with $|f(x) - f(y)| \leq C \beta (x, y)$, there is a {\it complex-valued} function in $Lip(X)$ with $F|_{S} = f$ and $|F(x) - F(y)| \leq 2C \beta (x, y)$ for all $x, y$ in $X$.
 
 \medskip
 
 {\bf Proof}.  We first check that the real and imaginary parts of $f$  are in $Lip(S)$ with the same Lipschitz constant $C$ as $f$.  By Proposition 1, there are real-valued $U,V$ in $Lip(X)$ with the same Lipschitz constant $C$ and such that $U|_{S} = Re(f), V|_{S} = Im(f)$.  Taking $F= U + iV$ gives the desired result.
 
 \medskip
 
Let $t_{0} = \sup \{ \beta(x, y): x, y \in X \}$.  If $\beta$ is unbounded, take $t_{0} = \infty$.  Assume that $t_{0} > 0$.  I can now prove the main result.
 
 \medskip
 
 \noindent
 {\bf Theorem 1}. On any metric space $(X, \beta)$, a complex-valued function $f$ is in $Lip_{c}(X)$ if and only if for every $\epsilon > 0$, there is a $ C = C(\epsilon) > 0$ so that
 
 $$
 |f(x) - f(y)| < \epsilon + C(\epsilon)\beta(x, y)
 \tag*
 $$
 
 \medskip
 
 \noindent
 for all $x, y$ in $X$.
 
 \medskip
 
 {\bf Proof.}  If $g$ is Lipschitz, with $|f(x) - g(x)| < \epsilon/2 $ for all $x$ in $X$ , then
 
 $$
 \align
 |f(x) - f(y)|&\leq |f(x) - g(x)| + |g(x) - g(y)| + |g(y) - f(y)|\\
 & < \epsilon + C\beta(x, y),
 \endalign
 $$
 
 \medskip
 
 \noindent
 where $C$ is a Lipschitz constant for $g$.

\medskip

For the converse, suppose $f$ satisfies (*) for every $\epsilon > 0$.  Without loss of generality, we may choose $C(\epsilon)$ in $(*)$ with $C(\epsilon) > \epsilon / t_{0} $ so there are $x_{1}, x_{2} \in X$ with $\beta (x_{1}, x_{2}) \geq \epsilon / C(\epsilon )$.   By Zorn's Lemma, with  $t = \epsilon/C(\epsilon)$ there is a maximal $t$-separated subset of $X$, $S = S_{\epsilon}$, which contains $x_{1}, x_{2}$.   For any $x, y$ in $S$ with $x \not= y$, we have $\beta (x, y) \geq t$ and by (*), 
 
 $$
 \align
{ {|f(x) - f(y)|}\over {\beta (x, y)}} &< {{\epsilon} \over { \beta(x, y)}} + C(\epsilon)\\
 & \leq \epsilon t^{-1} + C(\epsilon)\\
 &\leq 2 C(\epsilon),
 \endalign
 $$
 
 \medskip
 
 \noindent
 so $f|_{S} $ is Lipschitz with Lipschitz constant $2C(\epsilon)$.  
 
 \medskip
 
 \newpage
 
 \bigskip
 
 \bigskip
 
 By the Corollary to Proposition 1, $f|_{S}$ extends to a function $F$ which is in $Lip(X)$ and has Lipschitz constant $4 C(\epsilon)$.  For any $x$ in $X$, by maximality of $S$,  we may choose a $y$ in $S$ with $\beta (x, y) < t$.   Now $F(y) = f(y)$ so
 
 $$
 \align
 |F(x) - f(x)|& \leq |F(x) - F(y)| + |f(y) - f(x)|\\
 &\leq 4 C(\epsilon)\beta(x, y) + C(\epsilon)\beta(x, y) + \epsilon\\
 & < 6 \epsilon.
 \endalign
 $$

 \medskip
 
 \noindent
 Since $\epsilon > 0$ was arbitrary, the proof is complete.
 
 \medskip
 
 \noindent
 {\bf Corollary.} $Lip_{c}(X) \subset UC(X)$.
 
 \medskip
 
 {\bf Proof.}  Trivial.
 
 \medskip
 
 \noindent
 {\bf Remarks}.  The proof of Theorem 1 benefitted from a reading of [4, Proposition 2.1].  Functions in $Lip_{c}(X)$ can grow no faster than $f(x) = \beta (a, x)$ for any fixed $a$ in $X$.  In the next sections, I discuss some known examples  where $Lip_{c}(X) = UC(X)$.
 
 \bigskip
 
 {\bf 3. Complete Riemannian manifolds}.   I give a concise proof, using Theorem 1, of a known result [6, pp.286,289]. The key property of metric distance functions of complete (connected) Riemannian manifolds used here is:

\smallskip

\noindent
 geodesic completeness--between every two points $a, b$ there is a geodesic arc $\gamma$ of length $\beta (a, b)$.

 \medskip
 
  \noindent
 {\bf Proposition 2}. For any complete (connected) Riemannian manifold $(X, \beta)$, $Lip_{c}(X) = UC(X)$.
 
 \medskip
 
 {\bf Proof.}  For any $a \not= b$ in $X$, there is a geodesic segment $\gamma$ of length $\beta (a, b)$ joining $a$ to $b$.  For $f$ in $UC(X)$ and $\epsilon, \delta(\epsilon)$ as in the definition of uniform continuity above, let $N$ be the integer such that
 
 $$
 N \leq \beta(a, b) \delta(\epsilon)^{-1} < N + 1.
 $$
 
 \medskip
 
 \noindent
 Divide $\gamma$ into $N + 1$ equal-length segments, each of length less than $\delta(\epsilon)$. 
 
 \medskip
 
 \newpage
 
 \bigskip
 
 \bigskip
 
The triangle inequality then shows that
 
 $$
 \align
 |f(a) - f(b)| & < (N + 1) \epsilon\\ & \leq \beta(a, b) \delta(\epsilon)^{-1} \epsilon  + \epsilon \\
 & \leq C(\epsilon) \beta(a, b) + \epsilon,
 \endalign
 $$

\medskip

\noindent
where $C(\epsilon) = \epsilon \delta(\epsilon)^{-1}$.  Thus, (*) holds.  
 
\medskip

\noindent
{\bf Remark.}  This idea was used in [1, Lemma 2.1] when we considered the special case of {\it bounded symmetric domains} (BSD) $\Omega$ and obtained real-analytic Lipschitz approximants for all functions in $UC(\Omega)$.

\bigskip

 {\bf 4. Bounded symmetric domains}.  The bounded symmetric domains (BSD) in complex n-space ${\text{\bf C}}^{n}$ play a significant role in geometry and in representation theory [7].  These domains are all bounded open convex sets in ${\text{\bf C}}^{n}$ which carry intrinsic complete Riemannian (Bergman) metrics.  The prototype is just the hyperbolic metric on the open disc.  Using the results in [1], boundedness and compactness were determined for Toeplitz operators with uniformly continuous symbols on BSD's in [2].
 
 \medskip

There are two quite different natural metrics on BSD $\Omega$:  the Bergman metric, with distance function $\beta(\cdot, \cdot)$ and the restricted Euclidean metric from ${\text{\bf C}}^{n}$.
 The  two different corresponding notions of uniform continuity are related by the fact [5, p. 1167] that 
 $| x - y| \leq C_{\Omega} \beta(x, y)$ so that $UC(\Omega)_{|\cdot|} \subset UC(\Omega)_{\beta}$.  This provides a useful source of bounded functions in $UC(\Omega)_{\beta}$, which also includes unbounded functions like $f(z)= \beta(a, z)$ for any fixed $a$ in $\Omega$.  It follows easily from equation (*) of Theorem 1 that functions in $UC(\Omega)_{\beta}$ grow no faster than $\beta(a, z)$.  Finally, we observe that [3, Theorem E]  $\beta(a, z)$ is of slow growth near the boundary of $\Omega$: it is in $L^{p}(\Omega, dv)$ for {\it all} $p > 0$.

 \medskip
 
 \newpage
 
 \bigskip
 
 \bigskip
  
 \centerline{\bf References}
  
  \bigskip
  
 \noindent
[1] \ \ Bauer, W. and Coburn, L. A., Heat flow, weighted Bergman spaces, and real analytic Lipschitz approximation, {\it Journal fur die Reine und Angewandte Mathematik} (2015) 225-246.
 
 \medskip
 
\noindent
[2] \ \ Bauer, W. and Coburn, L.A.,  Toeplitz operators with uniformly continuous symbols, {\it Integral equations and operator theory} 83 (2015) 25-34.

\medskip

\noindent
[3] \ \ Bekolle, D., Berger, C. A., Coburn, L. A., Zhu, K. H., BMO in the Bergman metric on bounded symmetric domains, {\it Journal of functional analysis} 93 (1990)310-350.

\medskip

\noindent
[4] \ \ Benyamini, Y. and Lindenstrauss, J., {\it Geometric nonlinear functional analysis}, AMS Colloquium Publications 48 (2000) Providence, RI.

\medskip

\noindent
[5] \ \ Coburn, L. A., Sharp Berezin Lipschitz estimates, {\it Proceedings of the AMS} 135 (2007) 1163-1168.

\medskip

\noindent
[6] \ \ Garrido, M. I. and Jaramillo, J. A., Lipschitz-type functions on metric spaces, {\it J. Math. Analysis and Applications} 340 (2008) 282-290.

\medskip

\noindent
[7] \ \ Helgason, S., {\it Differential geometry, Lie groups, and symmetric spaces}, AMS Graduate Studies in Mathematics 34 (2001) Providence, RI.

\medskip

\noindent
[8] \ \ McShane, E. J., Extensions of range of functions, {\it Bulletin of the AMS} 40 (1934) 837-842.

\bigskip

\bigskip

\noindent
Department of Mathematics, SUNY at Buffalo, Buffalo, New York 14260, USA

\smallskip

\noindent
e-mail address: lcoburn$\@$buffalo.edu
\smallskip

\noindent
version: 5/13/2021

\vfill\eject\end